\documentclass[10pt,arxiv]{agn_article}
\usepackage{amsmath}
\usepackage{amsfonts}
\usepackage{amssymb}
\usepackage{subcaption}
\usepackage{graphics,epsfig}
\usepackage{color}
\usepackage[bibtex,natbib]{agn_all}
\begin{document}
\title{\vspace{-2cm}On the dispersion of waves for the  linear thermoelastic relaxed micromorphic model}
\date{\today}
\author {%
		Aarti Khurana\thanks{%
			Aarti Khurana,\quad Department of Mathematics, Panjab University, Chandigarh - 160014, India; email: aarti@pu.ac.in(AK)%
		}\quad and\quad %
	Suman Bala\thanks{%
		 Suman Bala,\quad Department of Mathematics, Panjab University, Chandigarh - 160014, India; 
		 email: sumanl$_{-}$75@pu.ac.in(SB)%
	}\quad and\quad %
	Hassam Khan\thanks{%
				corresponding author: Hassam Khan,\quad Lehrstuhl f\"{u}r Nichtlineare Analysis und Modellierung, Fakult\"{a}t f\"{u}r	Mathematik, Universit\"{a}t Duisburg-Essen, Thea-Leymann Str. 9, 45127 Essen, Germany, email: hassam.khan@stud.uni-due.de%
		}\\ \quad and\quad %
			Sushil K. Tomar\thanks{%
			 	Sushil K. Tomar,\quad Department of Mathematics, Panjab University, Chandigarh - 160014, India; email: sktomar@pu.ac.in(SKT)%
			}\quad and\quad %
				Patrizio Neff\thanks{%
		Patrizio Neff,\quad Head of Lehrstuhl f\"{u}r Nichtlineare Analysis und Modellierung, Fakult\"{a}t f\"{u}r	Mathematik, Universit\"{a}t Duisburg-Essen, Thea-Leymann Str. 9, 45127 Essen, Germany, email: patrizio.neff@uni-due.de%
	}
}
\maketitle
\begin{abstract}
\noindent We present the complete set of constitutive relations and field equations for the linear thermoelastic relaxed micromorphic continuum and investigate its variants for wave propagation. It is found that the additional thermal effects give rise to new waves and generate couplings with longitudinal waves which are not existing in the relaxed micromorphic continuum without thermal effects. However, transverse waves go un-affected by the thermal properties. The dispersion curves have been computed numerically for a particular model and compared with those presented by Madeo et al.\,\cite{madeo2015wave}.
\clearpage
\end{abstract}
\bigskip
\tableofcontents
\section{Introduction}
The classical micromorphic continuum theory developed by Eringen \cite{Eringen1970balance} and Mindlin \cite{mindlin1964micro} is perhaps the most generalized micro-continuum theory within the context of microstructural theories. In this theory, the deformation of the micromorphic continuum consists of macro-deformation and micro-deformation parts. The macro-deformation part is concerned with the deformation of the micromorphic solid in the sense of deformation of a classical continuum, while the micro-deformation part is the deformation lying beyond the scope of a classical deformation. This micro-deformation (or micro-distortion) yields additional degrees of freedom, which includes micro-strains, micro-stretches and micro-shear. Thus, each particle may have independent degrees of freedom describing micro-distortion, in addition to the macro-deformation degrees of freedom.\\
The origin of microstructural theories (micromorphic theories) goes back to the brothers Cosserat and Cosserat \cite{cosserat1909theorie}. In the Cosserat theory, each particle of the continua is taken to possess independent orientations along with macroscopic translation. The orientations describe the three degrees of freedom which are fully provided by three deformable directors associated with the center of mass of each infinitesimal particle of the continuum. Many theories have originated from the Cosserat theory,\, see e.g.\,\cite{,gunther1958statik,grioli1960elasticita,truesdell1960classical,toupin1962elastic,toupin1964theories,mindlin1964micro,jeong2010existence,neff2009new,jeong2009numerical}, Eringen and his co-worker\,\cite{eringen1964nonlinear,Eringen1970balance},\,Iesan \cite{iecsan2001theory,iecsan2002micromorphic} and Gale\c s\,\cite{galecs2010nonlinear} among others.\\
\indent  The development of the relaxed micromorphic model is a recent endeavor in view of a drastically simplified model of micromorphic type developed by Eringen. The classical micromorphic theory accounts for translational degrees of freedom at the macro level and micro-distortion degrees of freedom at the micro level. The force stress and the couple stress acting at the surface element of a classical micromorphic material is in general non-symmetric. However, the asymmetry of these tensors does not necessarily arise due to the microstructure of the body but comes into play due to some extra assumptions on the constitutive level. Neff et al.\,\cite{neff2014unifying} developed  their model by modifying these assumptions on the constitutive relations so that the force stress tensor may remain symmetric and the number of constitutive coefficients are reduced drastically, but the kinematic degrees of freedom stay the same. In addition, the curvature contribution, which is classically based on the full gradient $\nabla P$ of the microdistortion $P$ is "relaxed" to a dependence on only the  Curl $P$ of the microdistortion $P$. This specification completely changes the modeling capabilities  of the formulation. Indeed, arbitrary small samples in a classical micromorphic model respond with an infinite stiffness  against inhomogeneous deformation as  does a general gradient elasticity approach, while the relaxed micromorphic  model just turns into a particular stiffer classical  linear elastic  response. This new model is now called the `relaxed micromorphic continuum model' in the literature. Papers by Neff et al.\,\cite{neff2014unifying,madeo2016first,d2017panorama,agostino2019,neff2017real,neff2019identification} and Ghiba et al.\,\cite{ghiba2015relaxed}, contain the prominent foundation and explanations of the relaxed micromorphic continuum theory. They pointed out that in the new relaxed micromorphic model, the free energy is not uniformly point-wise positive-definite, but the new energy is positive semi-definite only. If the $\text{Curl}\,P$ contribution is neglected (zero characteristic length scale $\text{L}_c=0$)\,then the new model bears some resemblance to an early model of Hlav{\'a}{\'c}ek \cite{HLAVACEK19751137}. The extension to a transparent anisotropic formulation for the relaxed micromorphic is discussed in \cite{barbagallo2017transparent}. Madeo et al.\,\cite{madeo2015wave} have investigated the propagation of waves in an unbounded relaxed micromorphic continuum with 6 elastic parameters and explored the frequency range (band-gap) in which waves cannot be transmitted. These band-gaps appear due to the  presence of a coefficient occurring in the constitutive relations (the Cosserat couple modulus $\mu_c\geq 0$). This coefficient is also found to be responsible for the  non-symmetry of the force stress tensor, see \cite{eringen1969micromorphic} for a  related model. Later Madeo et al.\,\cite{madeo2016reflection} studied transmission and reflection of waves in band-gap mechanical metamaterials via the relaxed micromorphic model. Recently, Barbagallo et al.\,\cite{barbagallo2018relaxed} investigated the  transient wave propagation of the relaxed micrmorphic model in anisotropic band-gap metastructures. They showed that the relaxed micromorphic model gives a correct description of the pulse propagation in the frequency band-gap and at frequencies intersecting the optical branches, unlike Cauchy theory which only describes the  overall behavior of metastructures  at low frequencies. In addition,  Aivaliotis  et al.\,\cite{AIVALIOTIS201999}  considered the low-and high-frequency Stoneley waves, reflection and transmission at a Cauchy/relaxed micromorphic interface.\\
\indent \qquad In the relaxed micromorphic body, the kinematic variables are the classical displacement vector   ${\bf u}\!:\Omega\subset \mathbb{R}^3\rightarrow \mathbb{R}^3\notag $ and the micro-distortion tensor  $ P:\Omega\subset \mathbb{R}^3\rightarrow\mathbb{R}^{3\times3} $. The strain measures are taken to be 
\begin{align}\label{eqx01}
 e= \nabla {\bf u}-P,\qquad\textrm{sym}~P,\qquad\alpha=-\textrm{Curl}~P,
 \end{align}
 and the strain energy is given by
\begin{align}\label{eq01}
W&={}\frac{1}{2}\mathbb{C}^e_{klmn}(\textrm{sym}~e)_{kl}(\textrm{sym}~e)_{mn} +\frac{1}{2}\mathbb{C}^c_{klmn}(\textrm{skew}~e)_{kl}(\textrm{skew}~ e)_{mn}\nonumber\\&\quad +\frac{1}{2}\mathbb{C}^{\rm micro}_{klmn}
(\textrm{sym}~P)_{kl}(\textrm{sym}~P)_{mn} +\frac{1}{2} \mathbb{L}_{klmn}\alpha_{kl}\alpha_{mn}.
\end{align}
Here, the fourth-order elasticity tensors are defined as
\begin{align}
\mathbb{C}^e &:{}\mathrm{\Omega} \mapsto L(\mathrm{Sym}(3), \mathrm{Sym}(3)), \quad \mathbb{C}^c:\mathrm{\Omega} \mapsto L(\mathfrak{so}(3), \mathfrak{so}(3)),\nonumber\\\quad \mathbb{C}^{\rm micro}&:{}\mathrm{\Omega} \mapsto L(\mathrm{Sym}(3), \mathrm{Sym}(3)), \qquad~~ \mathbb{L}:\mathrm{\Omega} \mapsto L(\mathbb{R}^{3 \times 3}, \mathbb{R}^{3 \times 3}),
\end{align}
 where  $$\mathrm{Sym}(3):=\{X \in \mathbb{R}^{3 \times 3}| X^T=X\}, \qquad \mathfrak{so}(3):=\{X \in \mathbb{R}^{3 \times 3}| X^T=-X\}.$$
The non-symmetric parts of $P$ in the model are entirely due to moment stresses and applied moments.\\
Note that the skew-symmetric part of $P$ can not appear as individual contribution in \eqref{eq01} due to (linearized) frame indifference, i.e. the invariance
\begin{align}\label{eqo1}
\nabla {\bf u}\rightarrow\nabla{\bf u}+\overline{A},\qquad P \rightarrow P+\overline{A},
\end{align}
where $\overline{A}$ is any constant skew-symmetric matrix. We also note that $\nabla {\bf u}-P,\,\sym P$ and Curl $P$ are all invariant w.r.t. the transformation \eqref{eqo1}.\\\\
 \indent \qquad The aim of the present paper is concerned with the extension of the relaxed micromorphic model developed by Neff et al.\,\cite {neff2014unifying}\,to incorporate thermal effects. The relations and equations for the anisotropic relaxed micromorphic thermoelastic model have been taken from the Eringen-Mindlin model by considering the strain measures given in \eqref{eqx01} and considering the mixed terms with the  temperature variable only in the strain energy function. It is shown that the terms containing thermal effects are occurring in all the constitutive relations and the governing equations for the anisotropic micromorphic continuum. The constitutive relations and equations developed by Neff et al.\,\cite{neff2014unifying} are recovered in the absence of thermal effects, showing the validity of the derived equations. The proposed model may be useful for the description of metamaterials conducting heat.
\section{Conservation laws}
Consider a thermoelastic relaxed micromorphic continuum body occupying a finite region $\mathrm{\Omega},$ bounded by a piecewise smooth surface $\partial \mathrm{\Omega}$. The structure of the material is such that its particles can rotate, shear, stretch and shrink. With reference to a fixed system of rectangular coordinate axes $OX_{i} (i=1,2,3),$ the motion of the body is identified by the displacement $\mathbf{u}$ and micro-distortion $P$ of its macroscopic material points, defined by $\mathbf{u}=(u_{k}): \mathrm{\Omega}\times [0,t] \mapsto \mathbb{R}^{3}$ and $P = (P_{kl}): \mathrm{\Omega}\times [0,t] \mapsto \mathbb{R}^{3 \times 3},$ respectively. Thus, the quantities $(\mathbf{u},P)$ define the kinematic variables of the considered continuum model . We shall denote the components of the displacement vector by $u_k$ and $P$ is a non-symmetric second order micro-distortion tensor.
\noindent  We consider a set of constitutive variables for the relaxed micromorphic thermoelastic model, given by $$\{e, \textrm{sym}~P, \alpha, \theta \},$$
where
\begin{equation}\label{eq1}
e=\nabla \mathbf{u}-P, \quad \alpha = -\textrm{Curl}~P,
\quad \theta=\theta_{0}+\vartheta.\notag
\end{equation}
Here, $\theta$ is the temperature variable and $\vartheta$ is a small departure from the ambient temperature $\theta_{0}$.\\\\
For the relaxed micromorphic body, the expression of kinetic energy per unit mass, $\mathcal{I}$, in terms of generalized velocities is given as
\begin{equation}\label{eq11}
\mathcal{I}= \frac{1}{2}\dot{u}_{l}\dot{u}_{l}+\frac{1}{2}\zeta\dot{P}_{mk}\dot{P}_{mk},
\end{equation}
where superposed dot denotes the time derivative.
We define the total energy
\begin{equation}\label{eq12}
\mathcal{U}=\int_\Omega \left(\mathcal{I}+\mathcal{E}\right)dv,\nonumber
\end{equation}
per unit mass, where $\mathcal{E}$ is the total internal energy per unit mass. For the thermoelastic relaxed micromorphic continuum, the principle of balance of energy is written as
\begin{equation}\label{eq13}
\frac{d}{dt}\int_{ \mathrm{\Omega}}{\rho\,(\mathcal{E} + \mathcal{I})}\,dv = \int_{\partial \mathrm{\Omega}}(\sigma_{kl}\dot{u}_{l}+\epsilon_{kba}m_{lb}\dot{P}_{la}+q_{k})\, da_{k}+ \int_{\mathrm{\Omega}}\rho\,  (f_{l}\dot{u}_{l}+M_{lk}\dot{P}_{lk}+h) dv,
\end{equation}
where $da_{k}$ is the surface element, $dv$ is the volume element, $\sigma_{kl}$ is the force stress tensor, $m_{kl}$ is the non-symmetric  second order moment stress tensor\footnote{In the classical micromorphic model, $m$ would be a third order tensor.}, $h$ denotes the heat supply per unit mass, $q_{k}$ is  the heat vector, $f_l$ is the body force, $M_{kl}$ is the body moment and the symbol $\epsilon$ is the standard Levi-Civita tensor.\\ \indent
Following the procedure adopted by Eringen \cite{eringen2012microcontinuum}, the balance laws related to conservation of linear momentum and angular momentum are obtained as
\begin{align}\label{12}
\hat{f}_{l} &\equiv{} \sigma_{kl,k}+ \rho\,(f_{l}-\ddot{u}_{l})=0,\notag
\\ 
\hat{M}_{kl} &\equiv{} \epsilon_{lab}m_{kb,a} + \sigma_{lk}- s_{lk}+ \rho\,(M_{kl}- \Upsilon_{kl})=0,
\end{align}
where   $s_{kl}$ is any arbitrary symmetric tensor and   $\Upsilon_{kl}$ is the spin tensor given by $\Upsilon_{kl}\simeq \zeta\ddot{P}_{km}.$\\
\indent
In order to obtain the energy balance law, we apply the Transport and Green-Gauss theorems on equation \eqref{eq13} and utilizing equations \eqref{eq11} and \eqref{12}, we obtain
\begin{equation}\label{14}
-\rho \,\dot{\mathcal{E}}+ \sigma_{kl}(\dot{u}_{l,k}- \dot{P}_{lk})+ \epsilon_{lab}m_{kl}\dot{P}_{ka,b} +s_{kl}\dot{P}_{lk}+q_{k,k}+ \rho \,h= 0.\,
\end{equation}
Equations \eqref{12} and \eqref{14} constitute the local balance laws for the relaxed thermo-micromorphic medium. Equation \eqref{14} is the well known first law of thermodynamics.
\section{Constitutive relations}
\indent The second law of thermodynamics is given by
\begin{align}\label{slt}
\frac{d}{dt}\int_\Omega \rho\,\eta\, dv-\int_{\partial\Omega}\frac{1}{\theta}q_k\,n_k\,d{\bf a}-\int_{\Omega} \frac{\rho h}{\theta} dv  \geq 0,\,
\end{align}
where $\eta=\!_R \eta\, +_D \!\eta $ is the entropy density per unit mass; $_R \eta$ and  $_D\eta$ correspond to the static and dynamic parts of $\eta.$
Applying the Transport and Green-Gauss Theorems, we obtain the local form of the second law of thermodynamics as
\begin{eqnarray}\label{slt1}
\rho \,\dot\eta -\left(q_k/\theta\right)_{\!,k} -\rho\, h/\theta \geq 0.\,
\end{eqnarray}
Using the Helmoltz free energy $\psi$ given by
$\psi = \mathcal{E}-\theta \,\eta$ and eliminating $h$ from \eqref{14} and \eqref{slt1}, we obtain the Clausius-Duhem (C-D) inequality
\begin{align}\label{cdi}
 -\rho\,(\dot\psi+\eta\,\dot\theta)+\sigma_{kl}(\dot{u}_{l,k}- \dot{P}_{lk})+ \epsilon_{lab}\,m_{kl}\,\dot{P}_{ka,b} +s_{kl}\,\dot{P}_{lk}+\frac{1}{\theta}q_k \,\theta_{\!,k}\geq 0.\, 
 \end{align}
 Assuming that $\rho_0\,\psi=\Sigma(e,\mathrm{sym P}, \alpha, \theta_{\!,i} , \dot\theta,\theta; {\bf X}),$ the above inequality becomes
\begin{align}
\left(-\rho\,\eta-\frac{\rho}{\rho_0}\frac{\partial\Sigma}{ \partial \theta}\right)& \dot\theta+ \left(\sigma_{kl}-\frac{\rho}{\rho_0}\frac{\partial\Sigma}{\partial e_{kl}}\right)(\dot{u}_{l,k}- \dot{P}_{lk})+
\left(m_{kl}-\frac{\rho}{\rho_0}\frac{\partial\Sigma}{\partial \alpha_{kl}}\right)\epsilon_{lab}\dot{P}_{ka,b} &\qquad \nonumber \\
 +&\left(s_{kl}-\frac{\rho}{\rho_0}\frac{\partial\Sigma}{\partial(\mathrm{sym~P})_{kl}}\right)\dot{P}_{lk}
-\frac{\rho}{\rho_0}\frac{\partial\Sigma}{\partial \theta_{\!,k}}\dot\theta_{\!,k}-\frac{\rho}{\rho_0}\frac{\partial\Sigma}{\partial \dot\theta}\ddot\theta+\frac{1}{\theta}q_k\theta_{\!,k}\geq 0.
\end{align}
This inequality is linear in $\ddot \theta,~\dot{e}, \dot{P}, \dot{\alpha}$ and $\dot\theta_{\!,k}.$ It must remain positive for all independent variations of these quantities. But this is impossible unless
\begin{equation}\label{cr}
\left. \begin{aligned}
_R \eta & ={} -\frac{\rho}{\rho_0}\frac{\partial\Sigma}{\partial \theta},\qquad \qquad~ \sigma_{kl}=  \frac{\rho}{\rho_0}\frac{\partial\Sigma}{\partial e_{kl}}, \qquad \qquad \frac{\partial\Sigma}{\partial \theta_{,k}} = 0,\\ s_{kl}&={}\frac{\rho}{\rho_0}\frac{\partial\Sigma}{\partial(\mathrm{sym~P})_{kl}},\qquad  m_{kl} = \frac{\rho}{\rho_0}\frac{\partial\Sigma}{\partial \alpha_{kl}}, \qquad \qquad \frac{\partial\Sigma}{\partial \dot\theta} = 0,
\end{aligned}
\right\}
\end{equation}
and
\begin{align}\label{cdin}
-\rho _D \eta  \,\dot \theta +\frac{1}{\theta}q_k \,\theta_{,k}\geq 0.
\end{align}
 We see that $\Sigma$ is independent of the dynamic variables $\theta_{,k}$ and $\dot \theta$, thus $\psi$ is a static variable. Hence the only non-equilibrium part of the C-D inequality is given by \eqref{cdin}. Using  the Onsager postulate, we obtain
\begin{align}\label{qk}
q_k = \frac{\partial \Phi}{\partial (\theta_{\!,k}/\theta)},\qquad -\rho _D\eta = \frac{\partial \Phi}{\partial \dot\theta},
\end{align}
where $\Phi(=\Phi(\dot\theta,~ \theta_{,k}/\theta; (e, P, \alpha, \theta, {\bf X})))$ is the dissipation potential.\\
\indent Following Eringen \cite{eringen2012microcontinuum} and Neff et al. \cite{neff2014unifying}, the free energy $\psi$ and the dissipation potential $\Phi$ (with $_D \eta = 0$) for the thermoelastic  relaxed micromorphic continuum, are given by
\begin{align}\label{eq2}
\Sigma(e, P,\alpha, \theta) = \rho_{0}\psi = \frac{1}{2}\mathbb{C}^e_{klmn}(\textrm{sym}~e)_{kl}(\textrm{sym}~e)_{mn} +\frac{1}{2}\mathbb{C}^c_{klmn}(\textrm{skew}~e)_{kl}(\textrm{skew}~e)_{mn}\qquad&\nonumber \\+ \frac{1}{2}\mathbb{C}^{\rm micro}_{klmn}
(\textrm{sym}~P)_{kl}(\textrm{sym}~P)_{mn}+ \Sigma_{0}-\rho_{0}\,\eta_{0}\,\theta -\frac{\rho_{0}\,C_{0}}{2\theta_{0}}\theta^{2}\quad\notag\\- A_{kl}\theta (\textrm{sym}~e)_{kl} - G_{kl}\,\theta (\textrm{skew}~e)_{kl}- B_{kl} \,\theta
(\textrm{sym}~P)_{kl}\,\,\,&\\-D_{kl}\,\theta\,\alpha_{kl}+\frac{1}{2} \mathbb{L}_{klmn}\alpha_{kl}\,\alpha_{mn},\notag
\end{align}
and
\begin{align}\label{eq3}
\quad \Phi=\frac{1}{2\theta_{0}^{2}}\mathbb{K}_{kl}\,\theta_{\!,k}\theta_{\!,l}.
\end{align}
Some of the constitutive coefficients given in \eqref{eq2} and \eqref{eq3} follow the symmetries as
\begin{equation}\label{eq4}
\left.\begin{aligned}
\mathbb{C}^e_{klmn}&={}\mathbb{C}^e_{mnkl}=\mathbb{C}^e_{lkmn},\quad \mathbb{C}^c_{klmn}=\mathbb{C}^c_{mnkl}=-\mathbb{C}^c_{lkmn}, \quad A_{kl}=A_{lk}, \quad G_{kl}=-G_{lk},\\
\mathbb{C}^{\rm micro}_{klmn}&={} \mathbb{C}^{\rm micro}_{mnkl}=\mathbb{C}^{\rm micro}_{lkmn},\quad
\mathbb{L}_{klmn}= \mathbb{L}_{mnkl}, \quad B_{kl}= B_{lk}, \quad \mathbb{K}_{kl}= \mathbb{K}_{lk}.
\end{aligned}
\right\}
\end{equation}
It is noted that $\Phi$ is non-negative with an absolute minimum at $\theta_{\!,k}/\theta=0$ and $\dot\theta=0,$ i.e. 
\begin{align}
\Phi(0,~{\bf 0}; (e, P, \alpha, \theta, {\bf X}))=0.\,\notag
\end{align} 
Employing the expressions of $\Sigma$ and $\Phi$ in \eqref{cr} and \eqref{qk}, we obtain the constitutive relations as
\begin{align}
(\textrm{sym}~\sigma)_{kl}= {}& \frac{\partial \Sigma}{\partial (\textrm{sym}~e)_{kl}}= -A_{kl}\theta +\mathbb{C}^e_{klmn}(\textrm{sym}~e)_{mn},\notag\\ 
(\textrm{skew}~\sigma)_{kl}= {}&  \frac{\partial \Sigma}{\partial (\textrm{skew}~e)_{kl}}= -G_{kl}\theta +\mathbb{C}^c_{klmn}(\textrm{skew}~e)_{mn},\notag
\\
s_{kl}={} & \frac{\partial \Sigma}{\partial (\textrm{sym}~P)_{kl}}= -B_{kl}\theta+\mathbb{C}^{\rm micro}_{klmn}(\textrm{sym}~P)_{mn},\notag
\\ 
m_{kl}= {}&  \frac{\partial \Sigma}{\partial \alpha_{kl}}= -D_{kl}\theta+\mathbb{L}_{klmn}\alpha_{mn},\label{eq7}
\\ 
q_{k}= {} &\frac{\partial \Phi}{\partial (\theta_{\!,k}/\theta)}= \frac{1}{\theta_{0}}\mathbb{K}_{kl}\theta_{\!,l},
\notag\\ 
\eta={} &\frac{-1}{\rho_{0}}\frac{\partial \Sigma}{\partial \theta}= \eta_{0}+\frac{C_{0}}{\theta_{0}}\theta +\frac{1}{\rho_{0}} [A_{kl}(\textrm{sym}~e)_{kl}+ G_{kl}(\textrm{skew}~e)_{kl} \nonumber\\ &\qquad\qquad\qquad\qquad\qquad\quad+ B_{kl}(\textrm{sym}~P)_{kl}+D_{kl}\alpha_{kl}].\,\notag
\end{align}

Here, $\rho \simeq \rho_0$ is considered. In view of \eqref{cdin} and $\eqref{eq7}_5$ , it can be seen that $\mathbb{K}_{kl}$ is positive semi definite.
\section{Field equations}
\subsection{Nonlinear heat conduction}
Using \eqref{eq7} into equations \eqref{12}, and \eqref{14} along with the symmetry of the constitutive coefficients given in \eqref{eq4}, we obtain the field equations in index notation as
\begin{align}
(\mathbb{C}^e_{klmn}+\mathbb{C}^c_{klmn})(u_{n,mk}- P_{nm,k})- A_{kl}\theta_{\!,k}- G_{kl}\theta_{\!,k}+ \rho(f_{l}-\ddot{u}_{l})= 0,\notag
\\ \label{16}
(B_{kl} -A_{kl}-G_{kl})\theta- \epsilon_{lab}D_{kb}\theta_{\!,a}- \mathbb{L}_{kbmn}\epsilon_{lab}\epsilon_{nij}P_{mj,ia}+ (\mathbb{C}^e_{klmn}+\mathbb{C}^c_{lkmn})u_{n,m}\qquad\,\,\,& \nonumber \\-(\mathbb{C}^e_{klmn}+\mathbb{C}^c_{lkmn}+ \mathbb{C}^{\rm micro}_{klmn})P_{nm}+ \rho(M_{kl}- \zeta\ddot{P}_{kp})= 0,\\
-\frac{\rho_0\, C_0}{\theta_0}\theta\dot{\theta}+ \theta(A_{kl}-B_{kl}+G_{kl})\dot{P}_{lk}-\theta(A_{kl}+G_{kl}) \dot{u}_{l,k}+ \theta D_{kl}\epsilon_{lab}\dot{P}_{kb,a}+ \frac{1}{\theta _{0}}\mathbb{K}_{kl}\theta_{\!,kl}+ \rho h=0.\notag
\end{align}
The inherent rationale of the relaxed micromorphic model allows to \textbf{determine a priori} the appearing fourth order static elasticity tensors $\mathbb{C}^e, \mathbb{C}^{\rm micro}$ from a knowledge of the geometry on the micro-structural  level of the micromorphic material,  see \cite{agostino2019}. A similar identification for the classical Eringen-Mindlin micromorphic model is impossible. A further decisive advantage as compared to Eringen's approach is that the used elasticity tensors all have the classical  Voigt-format from linear elasticity, see \cite{barbagallo2017transparent}.

\subsection{Linearized heat conduction}
In terms of $\vartheta$ with $\theta=\theta_0+\vartheta $, equations \eqref{16} take the form
\begin{align}
(\mathbb{C}^e_{klmn}+\mathbb{C}^c_{klmn})(u_{n,mk}- P_{nm,k})- A_{kl}\theta_{\!,k}- G_{kl}\theta_{\!,k}+ \rho(f_{l}-\ddot{u}_{l})= &\,0,\notag\\
(B_{kl} -A_{lk}-G_{kl})(\theta_0+\vartheta)- \epsilon_{lab}D_{kb}\vartheta_{\!,a}- \mathbb{L}_{kbmn}\epsilon_{lab}\epsilon_{nij}P_{mj,ia}+ (\mathbb{C}^e_{klmn}+{C}^c_{lkmn})u_{n,m}\qquad&  \notag\\\mathbb-(\mathbb{C}^e_{klmn}+\mathbb{C}^c_{lkmn}+ \mathbb{C}^{\rm micro}_{klmn})P_{nm}+ \rho(M_{kl}-\zeta\ddot{P}_{kp})=&\, 0,\\
-\frac{\rho_0\, C_0}{\theta_0}(\theta_0+\vartheta)\dot{\vartheta}+ (\theta_0+\vartheta)(A_{kl}-B_{kl}+G_{kl})\dot{P}_{lk}-(\theta_0+\vartheta) (A_{kl}+G_{kl})\dot{u}_{l,k}\quad& \notag\\+(\theta_0+\vartheta) D_{kl}\epsilon_{lab}\dot{P}_{kb,a}+ \frac{1}{\theta _{0}}\mathbb{K}_{kl}\vartheta_{\!,kl}+ \rho\, h=&\,0.\notag
\end{align}
 In the following, we omit all the non-linear terms involving $\vartheta$, which lead to \\
 \begin{align}
 (\mathbb{C}^e_{klmn}+\mathbb{C}^c_{klmn})(u_{n,mk}- P_{nm,k})- A_{kl}\theta_{\!,k}+ \rho\,(f_{l}-\ddot{u}_{l})=& \,0,\notag
 \\ 
 (B_{kl} -A_{lk}-G_{kl})\,\theta- \epsilon_{lab}D_{kb}\vartheta_{\!,a}- \mathbb{L}_{kbmn}\epsilon_{lab}\epsilon_{nij}P_{mj,ia}+ (\mathbb{C}^e_{klmn}+\mathbb{C}^c_{lkmn})u_{n,m}\qquad& \nonumber \\-(\mathbb{C}^e_{klmn}+\mathbb{C}^c_{lkmn}+ \mathbb{C}^{\rm micro}_{klmn})P_{nm}+ \rho\,(M_{kl}-\zeta\ddot{P}_{kp})=& \,0,\\
 -\rho_0\, C_0\,\dot{\theta}+ \theta_0(A_{kl}-B_{kl}+G_{kl})\dot{P}_{lk}-\theta_0 A_{kl}\dot{u}_{l,k}-\theta_0 G_{kl}\dot{u}_{l,k}+ \theta_0 D_{kl}\epsilon_{lab}\dot{P}_{kb,a}\qquad&  \notag\\+ \frac{1}{\theta _{0}}\mathbb{K}_{kl}\theta_{\!,kl}+ \rho\,h=&\,0.\notag
 \end{align}
These are the local field equations for the linearized anisotropic  thermoelastic relaxed micromorphic solid. Now, we shall obtain the field equations for the isotropic solid with and without Cosserat couple modulus $\mu_c>0$ .
\subsection{Isotropic case}
\subsubsection{\kern-0.7emThe case for positive Cosserat couple modulus $\mu_c >0$}
\indent To write down the constitutive relations and field equations for the isotropic linear relaxed thermo-micromorphic continuum, we must take
\begin{align}
A_{kl}&= C_{1}\,\delta_{kl}, \quad B_{kl}= C_{2}\,\delta_{kl}, \quad D_{kl}= C_{3}\,\delta_{kl},  \quad \mathbb{K}_{kl}= C_4\,\delta_{kl},\notag \\
\mathbb{C}^e_{klmn}&=\mu_e(\delta_{km}\delta_{ln}+\delta_{kn}\delta_{lm})+\lambda_e\delta_{kl}\delta_{mn},\quad \mathbb{C}^e\sym X =2\mu_e\sym X+\lambda_e\tr(\sym X)\id ,\notag\\
\mathbb{C}^c_{klmn}&={} \mu_c(\delta_{km}\delta_{ln}-\delta_{kn}\delta_{lm}),\qquad \mathbb{C}^c\skew X=2\mu_c\skew X, \\ 
\mathbb{C}^{\rm micro}_{klmn}&=\mu_{\rm micro}(\delta_{km}\delta_{ln}+\delta_{kn}\delta_{lm}) +\lambda_{\rm micro}\delta_{kl}\delta_{mn},\quad \mathbb{C}^{\rm micro}\sym X =2\mu_{\rm micro}\sym X+\lambda_{\rm micro}\tr(\sym X)\id, \notag \\ 
\mathbb{L}_{klmn}&={} \overline{\alpha}_1\,\delta_{kl}\delta_{mn}+\overline{\alpha}_2\,\delta_{km}\delta_{ln}+\overline{\alpha}_3\,\delta_{kn}\delta_{lm}.\notag
\end{align}
The curvature energy expression takes the format 
\begin{align}
\frac{1}{2}\langle{\mathbb{L}\alpha,\alpha}\rangle=\frac{1}{2}\{\overline{\alpha}_1(\tr(\alpha))^2+\overline{\alpha}_2\vert\alpha\vert^2+ \overline{\alpha}_3\langle{\alpha,\alpha^T}\rangle\},\,
\end{align}
in which the parameters $\overline{\alpha}_1,\overline{\alpha}_2,\overline{\alpha}_3$ do not have themselves (but a combination thereof) a transparent physical interpretation \,(see \eqref{321}). Note that, in the isotropic case $G_{kl}\equiv0$ necessarily.\\
 Using these relations, the constitutive relations reduce to
\begin{align}
\sigma_{kl} &= (\mu_e+\mu_c)(u_{l,k}-P_{lk})+ (\mu_e-\mu_c)(u_{k,l}- P_{kl})+\lambda_e(u_{m,m}- P_{mm})\delta_{kl}- C_1 \theta \,\delta_{kl},\,\notag
\\
s_{kl}&= \mu_{\rm micro}(P_{lk}+ P_{kl})+\lambda_{\rm micro}P_{mm}\delta_{kl}- C_2 \,\theta \delta_{kl},\,\notag
\\\label{21}
m_{kl}&= -\overline{\alpha}_1\,\epsilon_{mab}P_{mb,a}\delta_{kl}- \overline{\alpha}_2\,\epsilon_{lab}P_{kb,a}-\overline{\alpha}_3\,\epsilon_{kab}P_{lb,a}- C_3 \,\theta \delta_{kl},\,
\\
\eta &={} \eta_{0}+ \frac{1}{\rho_0}\left[\frac{\rho_0\,C_{0}}{\theta_{0}}\theta + C_1u_{k,k}+ (C_2-C_1)P_{kk}- C_3\epsilon_{kab}P_{kb,a}\right],\,\notag
\\
q_k&={} \frac{C_4}{\theta_{0}}\theta_{\!,k}.\,\notag
\end{align}
 Plugging these expressions into \eqref{16}, the field equations for the isotropic relaxed thermo-micromorphic solid reduce to
\begin{align}
- C_1 \theta_{\!,l}+ (\mu_e+\mu_c)(u_{l,kk}-P_{lk,k})  +(\mu_e-\mu_c)(u_{k,lk}- P_{kl,k})+\lambda_e(u_{m,mk}- P_{mm,k})\delta_{kl}\notag\\+ \rho\,(f_l-\ddot{u}_l)&=  0,\notag
\\ \label{24}
(C_2-C_1)\theta \delta_{kl} -  C_3 \epsilon_{lik} \theta_{\!,i} - \overline{\alpha}_1\epsilon_{mab}\epsilon_{lik} P_{mb,ai} - \overline{\alpha}_2\epsilon_{lij}\epsilon_{jab}P_{kb,ai} - \overline{\alpha}_3\epsilon_{lij}\epsilon_{kab}P_{jb,ai} \qquad&\nonumber \\+  (\mu_e+\mu_c)(u_{k,l}- P_{kl})+(\mu_e-\mu_c)(u_{l,k}-P_{lk})+\lambda_e(u_{m,m}- P_{mm})\delta_{kl}\quad& \\-  \mu_{\rm micro}(P_{lk}+ P_{kl})-\lambda_{\rm micro}P_{mm}\delta_{kl} + \rho(M_{kl} - \zeta\ddot{P}_{kl})&=0,\notag
\\ 
\rho \,C_0\,\dot{\theta} - \frac{C_4}{\theta_{0}}\theta_{\!,kk} + \theta_0\{(C_2 - C_1)\dot{P}_{kk} + C_1 \dot{u}_{k,k} - C_3 \epsilon_{kab}\dot{P}_{kb,a}\} - \rho\, h&=0 ,\notag
\end{align}

\subsubsection{\kern-0.7emThe case for zero Cosserat couple modulus $\mu_c =0$}
In this case, the corresponding governing equation $\eqref{16}_2$ for the anisotropic thermoelastic relaxed micromorphic continuum reduces to
\begin{align}\label{26}
(B_{kl} -A_{kl})\theta- \epsilon_{lab}D_{kb}\theta_{\!,a}- \mathbb{L}_{kbmn}\epsilon_{lab}\epsilon_{nij}P_{mj,ia}+ \mathbb{C}^e_{klmn}u_{n,m}-(\mathbb{C}^e_{klnm}+ \mathbb{C}^{\rm micro}_{klmn})P_{nm}\quad& \nonumber \\+ \rho\,( M_{kl}- \zeta\ddot{P}_{kp})&= {}0,\,
\end{align}
while the other two equations given by  $\eqref{16}_1$ and  $\eqref{16}_3$ will remain unchanged. The constitutive relation $\eqref{21}_1$ will now take the form
\begin{eqnarray}\label{27}
\sigma_{kl}=  \mu_e(u_{l,k}+ u_{k,l}-P_{lk}- P_{kl})+\lambda_e(u_{m,m}- P_{mm})\delta_{kl}- C_{1}\theta \delta_{kl}\in \Sym\,(3),\,
\end{eqnarray}
while the other constitutive relations given in $\eqref{21}$ will remain unchanged. The corresponding governing equations $\eqref{16}$ become
\begin{alignat}{2}
&&-C_1\theta_{\!,l}+ \mu_e(u_{l,kk}-P_{lk,k}+ u_{k,lk}- P_{kl,k})+ \lambda_e(u_{m,mk}- P_{mm,k})\delta_{kl}+\rho\,(f_l-\ddot{u}_l)&= 0,\notag\\
\label{28}
&&(C_2-C_1)\theta\delta_{kl}- C_3 \epsilon_{lik} \theta_{\!,i} -\overline{\alpha}_1\epsilon_{mab}\epsilon_{lik} P_{mb,ai} - \overline{\alpha}_2\,\epsilon_{lij}\epsilon_{jab}P_{kb,ai}- \overline{\alpha}_3\,\epsilon_{lij}\epsilon_{kab}P_{jb,ai}\quad&  \notag \\ && +\mu_e(u_{l,k}-P_{lk}+ u_{k,l}- P_{kl})+\lambda_e(u_{m,m}- P_{mm})\delta_{kl}- \mu_{\rm micro}(P_{lk}+ P_{kl})\,\,\,&\\&&- \lambda_{\rm micro}P_{mm}\delta_{kl} + \rho\,(M_{kl} - \zeta\ddot{P}_{kl})&=0,\notag\\
&&\rho \,C_0\dot{\theta} - \frac{C_4}{\theta_{0}}\theta_{\!,kk} + \theta_0\{(C_2 - C_1)\dot{P}_{kk} + C_1 \dot{u}_{k,k} - C_3 \epsilon_{kab}\dot{P}_{kb,a}\} - \rho\,h&=0.\notag
\end{alignat}
\section{A special case: no thermal effects}

By neglecting the thermal effects from the present model, we are  left with the relaxed micromorphic continuum model. For this purpose, we  set $C_1=C_2=C_3=C_4=0.$ In order to reduce the constitutive relations and field equations for the relaxed micromorphic model with three parameters having positive Cosserat couple modulus $\mu_c > 0,$ and to achieve the governing equations of the model presented by Neff et al. \cite{neff2014unifying}, we make use of the notations given by 
\begin{align}\label{321}
\overline{\alpha}_1 &= a_3-\frac{a_1}{3}, \quad \overline{\alpha}_2= \frac{a_1+a_2}{2}, \quad \overline{\alpha}_3=  \frac{a_1-a_2}{2}.\end{align}
\qquad \qquad \qquad \qquad \qquad$\iff a_1=\overline\alpha_2+\overline\alpha_3,\qquad  a_2=\overline\alpha_2-\overline\alpha_3,\qquad  a_3=\frac{3\overline\alpha_1+\overline\alpha_2+\overline\alpha_3}{3}. $ \\

In terms of these relations the quadratic  curvature expression can equivalently be represented in irreducible components with physical meaning as
\begin{align}\label{x90}
\frac{1}{2}\langle{\mathbb{L}\alpha,\alpha}\rangle=\frac{1}{2}\{a_1\vert \dev\,\sym\, \alpha\vert^2+a_2\vert\,\skew\,\alpha\vert^2+a_3(\tr(\alpha))^2\}.\,
\end{align}
\\ With these, the constitutive relations $\eqref{21}_{1-3}$ become
\begin{align}
\sigma_{kl}&= (\mu_e+\mu_c)(u_{l,k}-P_{lk})+(\mu_e-\mu_c)(u_{k,l}-P_{kl})+{} \lambda_e(u_{m,m}- P_{mm})\delta_{kl},\notag
\\\label{33}
m_{kl}&={} a_1 \left[\frac{1}{3}\epsilon_{mab}P_{mb,a}\delta_{kl}-\frac{1}{2}\epsilon_{lab}P_{kb,a}- \frac{1}{2}\epsilon_{kab}P_{lb,a}\right] \\& \qquad+ \frac{a_2}{2}(\epsilon_{kab}P_{lb,a}-\epsilon_{lab}P_{kb,a})   -a_3\,\epsilon_{mab}P_{mb,a}\delta_{kl},\notag
\\ 
s_{kl}&={} 2\mu_{\rm micro}\frac{(P_{lk}+ P_{kl})}{2}+\lambda_{\rm micro}P_{mm}\delta_{kl},\notag
\end{align}
and the governing equations \eqref{28} reduce to
\begin{align}
(\mu_e+\mu_c)(u_{l,kk}-P_{lk,k})+ (\mu_e-\mu_c)(u_{k,lk}-P_{kl,k})+\lambda_e(u_{m,mk}- P_{mm,k})\delta_{kl} \qquad\notag&\\+\rho\,(f_l-\ddot{u}_l)&=\,0,\notag
\\ \label{38}
-a_1\, \epsilon_{lij}\left(\frac{1}{2}\epsilon_{kab} P_{jb,ai}+\frac{1}{2}\epsilon_{jab}P_{kb,ai} - \frac{1}{3}\epsilon_{mab}P_{mb,ai}\delta_{kj}\right)- a_3\, \epsilon_{lik} \epsilon_{mab}P_{mb,ai} \qquad\notag&\\+\frac{a_2}{2} \epsilon_{lij}\left(\epsilon_{kab}P_{jb,ai}-\epsilon_{jab} P_{kb,ai}\right)+(\mu_e+\mu_c)(u_{k,l}-P_{kl})\quad&\\+(\mu_e-\mu_c)(u_{l,k}-P_{lk})\lambda_eu_{m,m}\delta_{kl}-\mu_{\rm micro}(P_{kl}+P_{lk})\,\,\,\notag&\\-(\lambda_e+\lambda_{\rm micro}) P_{mm} \delta_{kl}+ \rho\,(M_{kl} - \zeta\ddot{P}_{kl})&= 0.\notag
\end{align}

Note that the equation $\eqref{28}_3$ is automatically satisfied in the absence of thermal effects. The equations \eqref{33} are the constitutive relations and field equations for the relaxed micromorphic continuum model having three parameters, namely, $a_1$, $a_2$ and $a_3$ with positive Cosserat couple modulus $\mu_c > 0.$ \,\\
\indent Further, by setting the Cosserat couple modulus $\mu_c=0$ in \eqref{33} one can recover the constitutive relations and governing equations given by (2.24) and (2.25) earlier derived by Neff et al. \cite{neff2014unifying} for the relaxed micromorphic model with zero Cosserat couple modulus $\mu_c=0.$\,
\section{Plane wave propagation}
Here, we shall explore the possibility of one-dimensional wave propagation in the linear relaxed thermoelastic micromorphic solid of infinite extent. For this purpose, let us consider the propagation of a wave along the $x$-direction and that all the kinematical variables depend on the spatial coordinate $x$ only. As in Madeo et al.\,\cite{madeo2015band}, the component $P_{11}$ of the tensor $P$ can be written as
\begin{equation}
P_{11}=P^D+P^S,
\end{equation}
where
\begin{equation}\label{39}
P^S= \frac{1}{3}(P_{11}+P_{22}+P_{33}), \quad  P^D= (\textrm{dev}~ \textrm{sym} ~P )_{11}.
\end{equation}
We introduce the following notations
\begin{equation}\label{40}
P^V=P_{22}-P_{33},\quad \zeta_0=\ \rho\, \zeta, \quad (\textrm{sym}~P)_{1 \xi}=P_{(1 \xi)},\quad (\textrm{skew} P)_{1\xi}=P_{[1\xi]}, \quad \xi=1,2.
\end{equation}
Using the relations given in \eqref{39} and \eqref{40}, the field equations \eqref{24} give rise to the following sets of equations
\\ \noindent {\bf(i) \em  Coupled longitudinal waves : } 
\begin{equation}\label{42'}
\left.\begin{aligned}
\ddot{u}_1&={}\frac{-C_1}{\rho}\theta_{,1}+c^2_p u_{1,11}-\frac{2 \mu_e}{\rho}P^D_{,1}-\frac{3\lambda_e+2\mu_e}{\rho}P^S_{,1},\\
\ddot{P}^D&={}\frac{4\mu_e}{3\zeta_0}u_{1,1}+\frac{\overline{\alpha}_2-\overline{\alpha}_3}{3\zeta_0}P^D_{,11}-\frac{2(\overline{\alpha}_2-\overline{\alpha}_3)}{3\zeta_0}P^S_{,11}-w^2_s P^D,\\
\ddot{P}^S&={}\frac{(C_2-C_1)}{\zeta_0} \theta+\frac{3\lambda_e+2\mu_e}{3\zeta_0}u_{1,1}-\frac{\overline{\alpha}_2-\overline{\alpha}_3}{3\zeta_0}P^D_{,11}+\frac{2(\overline{\alpha}_2-\overline{\alpha}_3)}{3\zeta_0}P^S_{,11}-w^2_p P^S,\\
\ddot{P}_{[23]}&={}\frac{-C_3}{\zeta_0}\theta_{,1}+\frac{\overline{\alpha}_2+\overline{\alpha}_3+2\overline{\alpha}_1}{\zeta_0}P_{[23],11}-\frac{2\mu_c}{\zeta_0}P_{[23]},\\
\frac{C_4}{\theta_0}\theta_{,11}&={}\rho\, C_0\dot{\theta}+\theta_0 \left[3(C_2-C_1)\dot{P}^S+C_1\dot{u}_{1,1}+2C_3\dot{P}_{[23],1}\right],
\end{aligned}
\right\}
\end{equation}
\begin{align}
c^2_p={\frac{\lambda_e+2\mu_e}{\rho}},~~ \omega^2_s={\frac{2\,(\mu_e+\mu_{\rm micro})}{\zeta_0}},~~\omega^2_p=\frac{2\,(\mu_e+\mu_{\rm micro})+3\,(\lambda_e+\lambda_{\rm micro})}{\zeta_0}.\notag\,
\end{align}
Equations in \eqref{42'} are coupled in $u_1,~P^S,~P^D,P_{[23]}$ and $\theta$ in $x$-direction.\\
[1cm]{\bf(ii) \em  Coupled transverse waves:}
\begin{equation}\label{43a}
\left.\begin{aligned}
\ddot{u}_\xi&={}c^2_s u_{\xi,11}-\frac{2\mu_e}{\rho} P_{(1\xi),1}+\frac{2\mu_c}{\rho} P_{[1\xi],1}, \quad c^2_s={\frac{\mu_e+\mu_c}{\rho}}
\\
\ddot{P}_{(1\xi)}&={}\frac{\mu_e}{\zeta_0}u_{\xi,1}-w^2_s P_{(1\xi)}+\frac{\overline{\alpha}_2}{2\zeta_0}P_{(1\xi),11}+\frac{\overline{\alpha}_2}{2\zeta_0}P_{[1\xi],11},\\
\ddot{P}_{[1\xi]}&={}-\frac{\mu_c}{\zeta_0}u_{\xi,1}-\frac{2\mu_c}{\zeta_0}P_{[1\xi]}+\frac{\overline{\alpha}_2}{2\zeta_0}P_{(1\xi),11}+\frac{\overline{\alpha}_2}{2\zeta_0}P_{[1\xi],11}.
\end{aligned}
\right\}
\end{equation}
These are two sets of three coupled equations in $u_{\xi},~{P}_{(1\xi)}$ and ${P}_{[1\xi]}$ in the $\xi^{th}$ direction with $\xi=2,3.$\\[1cm]
\noindent {\bf(iii) \em  Uncoupled longitudinal waves:}
\begin{equation}\label{44}
\left.\begin{aligned}
\ddot{P}_{(23)}&={}c^2_m P_{(23),11}-w^2_s P_{(23)},
\\
\ddot{P}^V&={}c^2_m P^V_{,11}-w^2_s P^V, \quad c^2_m={\frac{\overline{\alpha}_2+\overline{\alpha}_3}{\zeta_0}}.
\end{aligned}
\right\}
\end{equation}
These equations are uncoupled in $P_{(23)}$ and $P^V.$\\
\indent For a time harmonic plane progressive wave in the positive direction of the  $x$-axis, we take the following form of various entities
\begin{equation}\label{46}
\left.\begin{aligned}
(u_1,~P^D,~P^S,~P_{[23]},~\theta)&={} Re\{(\beta_1,~\beta_2,~\beta_3,~\beta_4,~\beta_5) e^{\dot{\iota}(k x-\omega t)}\}, \\
(u_\xi,~P_{(1\xi)},~P_{[1\xi]})&={} Re\{(\gamma^\xi_{1},~\gamma^\xi_2,~\gamma^\xi_3) e^{\dot{\iota}(k x-\omega t)}\},~~\xi= 2,3,\\
(P_{(23)}, P^V)&={} Re\{(\beta_{(23)}, \beta^V)e^{\dot{\iota}(k x-\omega t)}\},
\end{aligned}
\right\}
\end{equation}
where $\beta_j~(j=1,...,5),$ $\gamma^\xi_l~(k=1,2,3),$ $\beta_{(23)}$ and $\beta^V$ are the amplitudes, $k$ is the wavenumber, $\omega$ is the angular frequency and $\dot{\iota} = \sqrt{-1}$ is the imaginary unit.\\
Now plugging the relevant quantities from \eqref{46} into the sets of equations given by \eqref{42'} and \eqref{43a}, one can obtain
\begin{equation}\label{47}
\textbf{A}_1 \cdot \boldsymbol{\beta}=0, \qquad \textbf{A}_\xi \cdot \boldsymbol{\gamma}^\xi=0, \quad \xi= 2,3,
\end{equation}
where $\textbf{A}_1= [a_{ij}]$ is a $5 \times 5$ matrix, whose non-zero elements are given by
\begin{align}
a_{11}&={}\omega^2-c^2_p k^2, \quad a_{12}= -\frac{2\dot{\iota} k \mu_e}{\rho},\quad a_{13}= -\frac{\dot{\iota} k(3\lambda_e+2\mu_e)}{\rho},\quad a_{15}= -\frac{\dot{\iota} k C_1}{\rho},\nonumber \\
a_{21}&={} \frac{4 \dot{\iota}k \mu_e}{3\zeta_0}, \quad a_{22}= \omega^2-\frac{k^2 (\overline{\alpha}_2-\overline{\alpha}_3)}{3\zeta_0}- \omega^2_s, \quad  a_{23}=\frac{2 k^2 (\overline{\alpha}_2-\overline{\alpha}_3)}{3\zeta_0},\nonumber\\
a_{31}&={}\frac{\dot{\iota} k(3\lambda_e+2\mu_e)}{3\zeta_0},\quad a_{32}=\frac{k^2(\overline{\alpha}_2-\overline{\alpha}_3)}{3\zeta_0}, \quad a_{33}= \omega^2-\omega^2_p-\frac{2k^2(\overline{\alpha}_2-\overline{\alpha}_3)}{3\zeta_0},\\ a_{35}&={}\frac{C_2-C_1}{\zeta_0},\quad
a_{44}=\omega^2- \frac{(\overline{\alpha}_2+\overline{\alpha}_3+2\overline{\alpha}_1)k^2}{\zeta_0}-\frac{2\mu_c}{\zeta_0}, \quad a_{45}= -\frac{\dot{\iota} k C_3}{\zeta_0}, \nonumber\\
a_{51}&={}k \theta_0 C_1 \omega,\quad a_{53}= -3 \dot{\iota} \omega(C_2-C_1) \theta_0, \quad a_{54}= 2 k C_3 \theta_0 \omega, \quad a_{55}= -\rho \,\dot{\iota} C_0 \omega+\frac{k^2C_4}{\theta_0},\nonumber
\end{align}
while the matrices $\textbf{A}_\xi$ are given by
\begin{align} 
\textbf{A}_2= \textbf{A}_3=\left({\begin{array}{ccc}
-\omega^2+ c_s^2 k^2 & \displaystyle\frac{2 \dot{\iota} k \mu_e}{\rho} & - \displaystyle\frac{2 \dot{\iota}k \mu_c}{\rho}\\~\\ -\displaystyle\frac{2\dot{\iota} k \mu_e}{\zeta_0} & -2 \omega^2+\displaystyle\frac{k^2 \overline{\alpha}_2}{\zeta_0}+ 2\omega^2_s & \displaystyle\frac{k^2 \overline{\alpha}_2}{\zeta_0}\\~\\\displaystyle\frac{2\dot{\iota} k \mu_c}{\zeta_0} & \displaystyle\frac{k^2 \overline{\alpha}_2}{\zeta_0} & -2\omega^2+\displaystyle\frac{k^2\overline{\alpha}_2}{\zeta_0}+ \frac{4\mu_c}{\zeta_0}
	\end{array}}\right),
\end{align}
and $${\boldsymbol \beta}=(\beta_1, \beta_2, \beta_3, \beta_4, \beta_5)^T, \quad {\boldsymbol \gamma}^\xi=(\gamma^\xi_1, \gamma^\xi_2, \gamma^\xi_3)^T.$$
For a non-trivial solution of the algebraic system \eqref{47}, we must have
\begin{equation}\label{48}
\mbox{det} \textbf{A}_{1}= 0,\quad \mbox{det} \textbf{A}_{2}= 0,\quad \mbox{det} \textbf{A}_{3}= 0.
\end{equation}
Similarly, inserting the relevant entities from \eqref{46} into equations in \eqref{44}, we obtain the following same dispersion relation given by \begin{equation}\label{49}
\omega=\sqrt{c^2_m k^2+\omega^2_s}.
\end{equation}
Equations  \eqref{48} and \eqref{49} are the dispersion relations of various existing waves in the  relaxed micromorphic model with three curvature parameters and Cosserat couple modulus $\mu_c \geq 0.$
\section{Dispersion relations under particular cases}
	1. In order to reduce the dispersion relations for the relaxed thermo-micromorphic model having three-parameters without Cosserat couple modulus $\mu_c,$ we
	substitute  $\mu_c=0$ in the elements of the equations given in \eqref{48}. Note that the dispersion relation \eqref{49} will remain unchanged. The dispersion relations is reduced from \eqref{48} together with \eqref{49} would represent those for this model.\\~\\
	2. To reduce the dispersion relations for the relaxed thermo-micromorphic model having one-parameter with and without Cosserat couple modulus $\mu_c,$ with reference to \eqref{x90} we consider the following three cases:\\\\
   (i) $a_1=a_2>0,\quad  a_3=\frac{a_1}{3},$\\
   (ii) $a_1>0, a_2=0,\quad a_3=\frac{a_1}{3},$\\
   (iii) $a_1=a_2 >0,\quad  a_3=0.$\\
\\\noindent	Within this setting, one can write the corresponding dispersion relations from \eqref{48} and \eqref{49} with and without Cosserat couple modulus $\mu_c$.
\section{Numerical results and discussion}
	To study thermal effects on dispersion curves in greater details, we have performed numerical computations for a specific model. The relevant values of different elastic moduli have been borrowed from Madeo et al.\,\cite{madeo2015wave} and those for the thermal parameters have been taken appropriately. These are given in the following table. \, \begin{center}
		\begin{tabular}{|c|c|c|}
			\hline
			Parameter & Value & Unit  \\ \hline
			$\lambda_e$ & 400 & MPa \\
			$\mu_e$ & 200 & MPa \\
			$\lambda_{\rm micro}$ & 100 & MPa \\
			$\mu_{\rm micro}$ & 100 & MPa \\
			$\mu_c$ & 440 & MPa \\
			$\rho$ & 2000 & Kg/m$^3$ \\
			$\zeta_0$ & $0.01$ & Kg/m \\
			$C_0$ & 206 & Cal~Kg$^{-1}$~$^o C^{-1}$ \\
			$C_1$ & 84 $\times 10^{3}$ & Cal~m$^{-3}$~$^o C^{-1}$ \\
			$C_2$ & 95 $\times 10^{3}$ & Cal~m$^{-3}$~$^o C^{-1}$ \\
			$C_3$ & 152  & Cal~m$^{-2}$~$^o C^{-1}$ \\
			$C_4$ & 16  & Cal~m$^{-1}~s^{-1}$ \\
			$\theta_0$ & 20 & $^o C$ \\
			$a_1$ & 2$\times 10^{-3}$ & MPa m$^2$\\
			$a_2$ & 2$\times 10^{-3}$ & MPa m$^2$ \\
			$a_3$ & 2$\times 10^{-3}$ & MPa m$^2$\\
			\hline
		\end{tabular}
	\end{center}
	With this numerical data, the dispersion equations \eqref{48} and \eqref{49} are solved through MatLab-Software. We have considered the following six thermoelastic relaxed micromorphic models:\\ \noindent  {\bf Model-I:} Thermoelastic relaxed micromorphic model having $a_1=a_2>0,$ $a_3=\displaystyle \frac{a_1}{3}$ with positive Cosserat couple modulus $\mu_c >0$,\\ \noindent  {\bf Model-II:} Thermoelastic relaxed micromorphic model having $a_1=a_2>0,$ $a_3=\displaystyle \frac{a_1}{3}$ with zero Cosserat couple modulus $\mu_c =0$,\\ \noindent {\bf Model-III:} Thermoelastic relaxed micromorphic model having $a_1>0, a_2=0,$ $a_3=\displaystyle \frac{a_1}{3}$ with positiveCosserat couple modulus $\mu_c >0$,\\ \noindent {\bf Model-IV:} Thermoelastic relaxed micromorphic model having $a_1>0, a_2=0,$ $a_3=\displaystyle \frac{a_1}{3}$ with zero Cosserat couple modulus $\mu_c =0,$\\\noindent  {\bf Model-V:} Thermoelastic relaxed micromorphic model having $a_1=a_2 >0,$ $a_3=0$ with positive Cosserat couple modulus $\mu_c >0$,\\ \noindent {\bf Model-VI:} Thermoelastic relaxed micromorphic model having $a_1=a_2 >0,$ $a_3=0$ with zero Cosserat couple modulus $\mu_c =0.$\\ The various dispersion curves obtained have been depicted graphically through Figures 1-7.\\
	\begin{figure*}
		\begin{center}
			\includegraphics[width=0.85\textwidth]{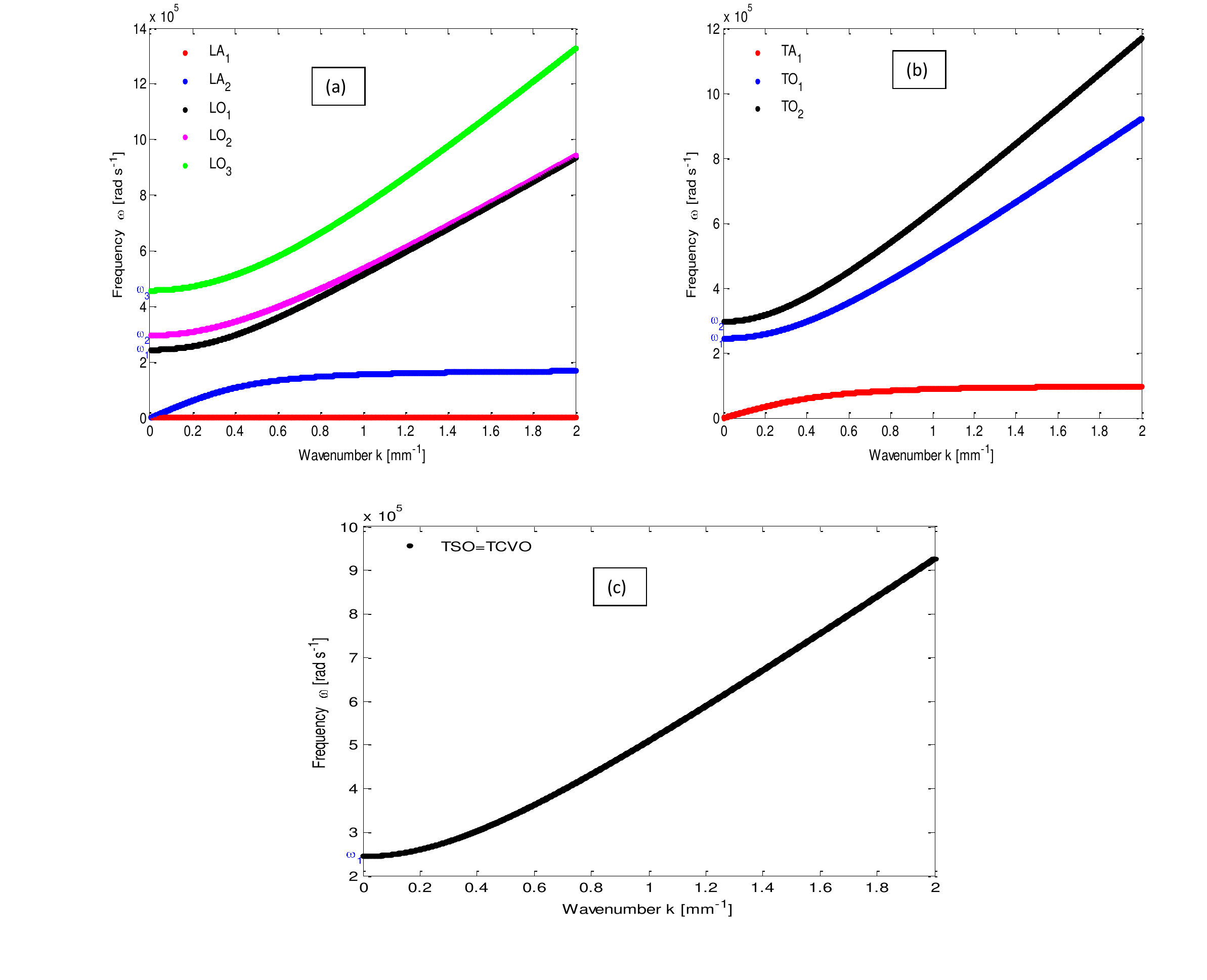}
			\caption{Dispersion curves for the thermoelastic relaxed micro-morphic model having $a_1=a_2>0,$ $a_3=\displaystyle \frac{a_1}{3}$ with positive Cosserat couple modulus $\mu_c > 0$. [(a) Longitudinal branches, (b) Transverse branches, (c) Uncoupled Optic
				branch]}
		\end{center}
	\end{figure*}
	\indent Figure 1(a) depicts the dispersion curves corresponding to longitudinal acoustic (LA) and longitudinal optic (LO) modes, while Figure 1(b) depicts that of corresponding to transverse acoustic (TA) and transverse optic (TO) modes. Figure 1(c) represents the dispersion curve of uncoupled optic modes in the relaxed thermo-micromorphic model-I with positive Cosserat couple modulus $\mu_c > 0.$
	\begin{figure*}
		\begin{center}
			\includegraphics[width=0.75\textwidth]{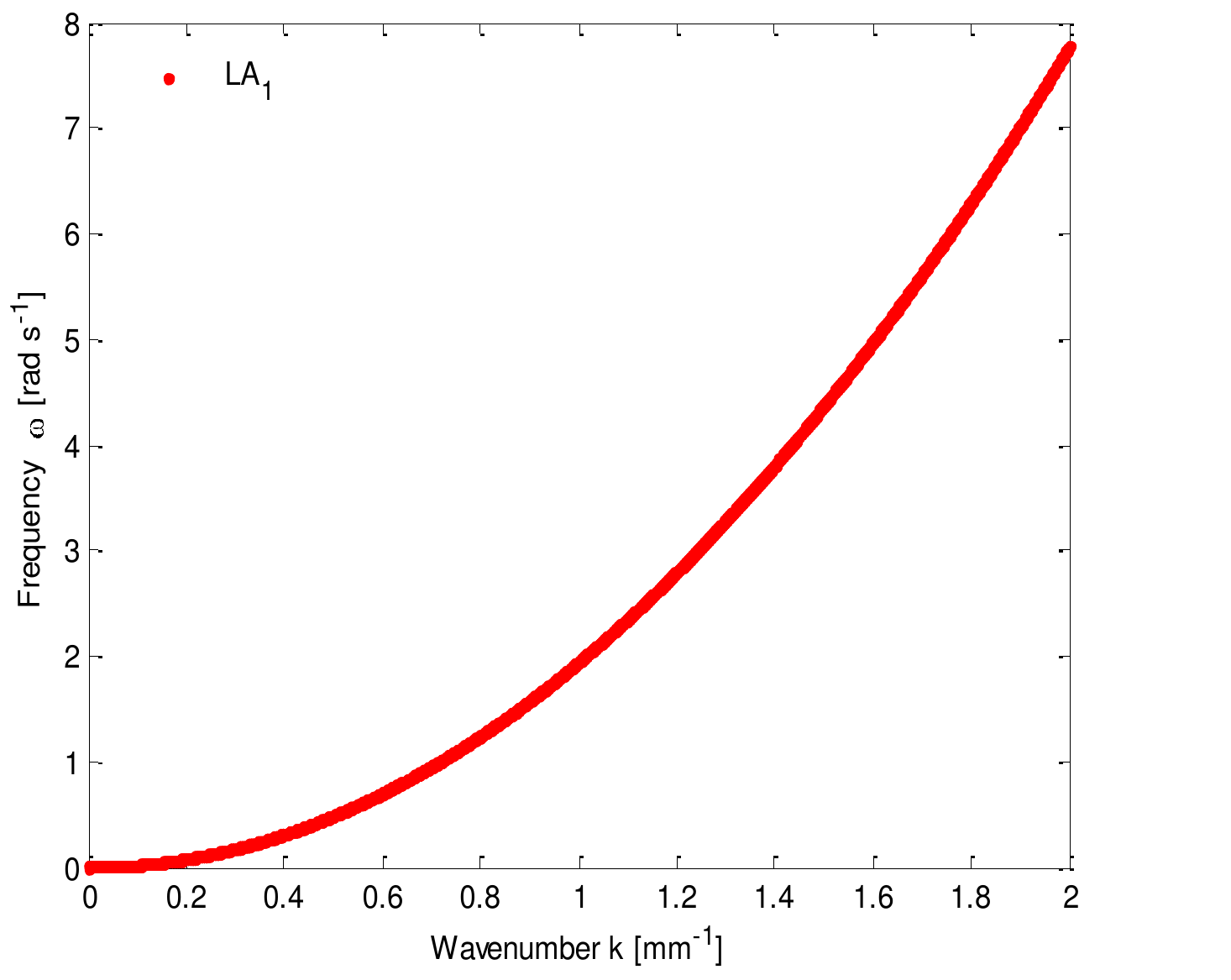}
			\caption{Zooming of first LA branch}
		\end{center}
	\end{figure*}
	Figure 2 depicts the zooming view of the first LA branch. It can be observed from Figure 1(a) that there exist two LA branches and three LO branches. One of the two LA branches is less dispersive than the other branch. The graph of this LA branch has been depicted separately through Figure 2, which shows that it is dispersive but less dispersive relative to the other LA branch. We also note that all the three LO branches are significantly dispersive. The cut-off frequencies of LO$_1,$ LO$_2$ and LO$_3$ branches are $2.4495 \times 10^5~rad/s,$ $2.9665 \times 10^5~rad/s$ and $4.5826 \times 10^5~rad/s,$ respectively. From Figure 1(b), we note that there is only one TA branch and two TO branches having cut-off frequencies, namely, $2.4495 \times 10^5~rad/s$ and $2.9665 \times 10^5~rad/s$. This model can be compared with that of Madeo et al.\,\cite{madeo2015wave} in contrast with thermal effects. First and foremost, it can be seen that the presence of thermal field in the model of Madeo et al.\,\cite{madeo2015wave} gives rise to (i) a new dispersion curve and (ii) transforms an uncoupled mode to a coupled longitudinal mode. The new dispersion curve corresponds to an LA mode and the uncoupled TRA branch of Madeo et al.\,\cite{madeo2015wave} transforms into a coupled first LO branch. There remains only one uncoupled TSO-TCVO branch. Note that the frequency band gap which arose in Madeo et al.\,\cite{madeo2015wave} is not at all influenced by the thermal field. Figure 3 depicts the dispersion curves of various existing waves in Model-II. The effect of the Cosserat couple modulus $\mu_c$ on the dispersion  curves can be noticed by comparing the respective graphs of Figures 1 and 3. On comparing these figures, one can notice that the presence of the Cosserat couple modulus $\mu_c$ is solely responsible for frequency band gaps. On comparing Figure 1(a) with Figure 3(a), it is seen that there exist three-LO and two-LA branches for positive Cosserat couple modulus $\mu_c >0$, while there exist two-LO and three-LA branches for zero Cosserat couple modulus $\mu_c=0.$ A similar phenomena is seen while comparing Figure 1(b) with Figure 3(b). On comparing Figure 1(c) with Figure 3(c), we see that there are no effects  of the Cosserat couple modulus   $\mu_c$ on TSO-TCVO branches same expression of the frequency given through \eqref{49} is independent of the Cosserat couple modulus $\mu_c.$\\
	\begin{figure*}
		\begin{center}
			\includegraphics[width=0.85\textwidth]{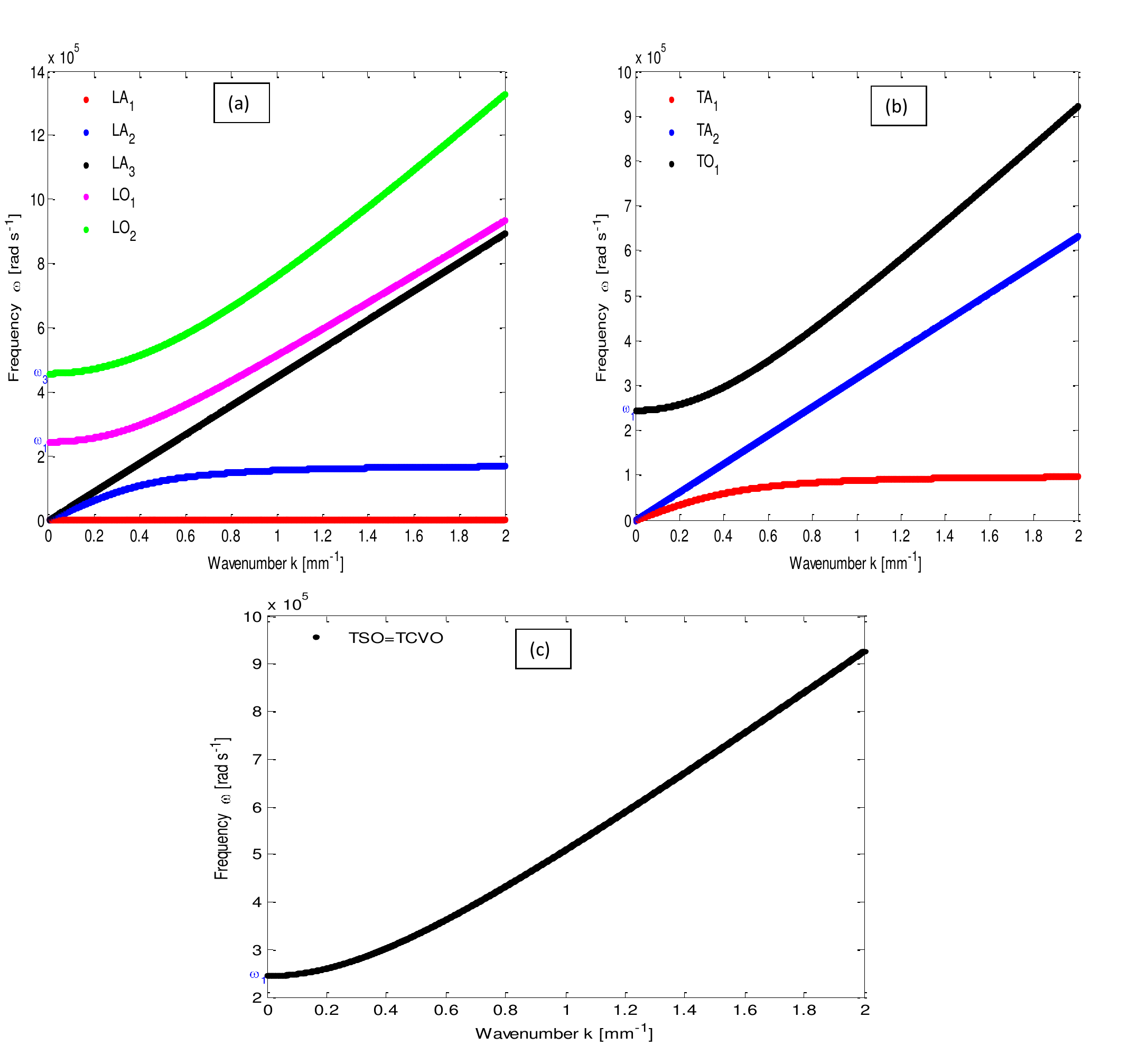}
			\caption{Dispersion curves in the thermoelastic relaxed micro-morphic model having $a_1=a_2>0,$ $a_3=\displaystyle \frac{a_1}{3}$ with zero Cosserat couple modulus $\mu_c = 0$. [(a) Longitudinal branches, (b) Transverse branches, (c) Uncoupled Optic
				branch]}
		\end{center}
	\end{figure*}
	\begin{figure*}
		\begin{center}
			\includegraphics[width=0.85\textwidth]{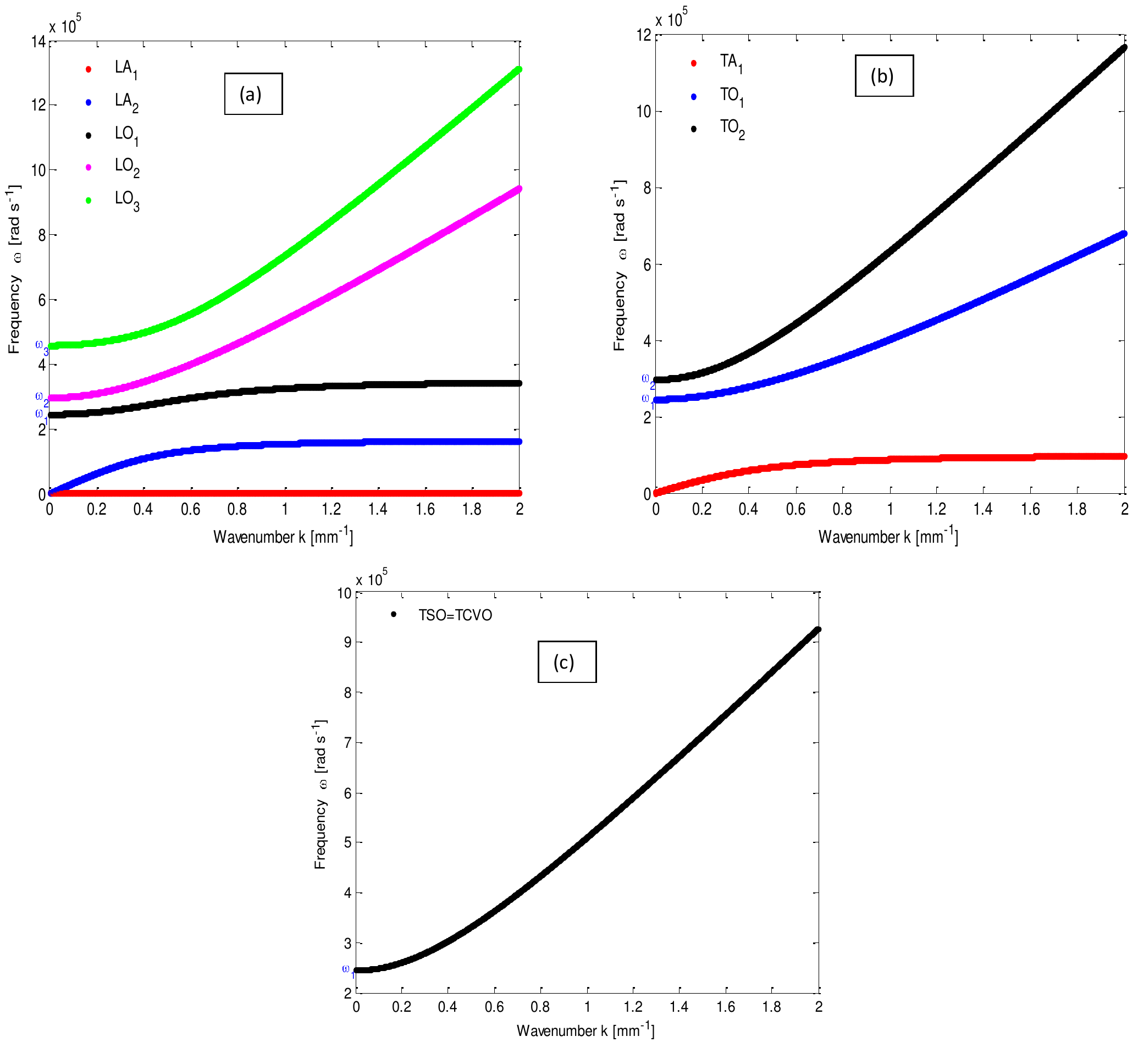}
			\caption{Dispersion curves in the thermoelastic relaxed micro-morphic model having $a_1>0, a_2=0,$ $a_3=\displaystyle \frac{a_1}{3}$ with positive Cosserat couple modulus $\mu_c > 0$. [(a) Longitudinal branches, (b) Transverse branches, (c) Uncoupled Optic
				branch]}
		\end{center}
	\end{figure*}
	
	\indent Figures 4 and 5 depict the dispersion curves of various acoustic and optic modes respectively, for Models-III and IV. From Figure 4, it is noted that there exists a band gap in Model-III, while no band gap is observed in Model-IV.\\
	\begin{figure*}
		\begin{center}
			\includegraphics[width=0.85\textwidth]{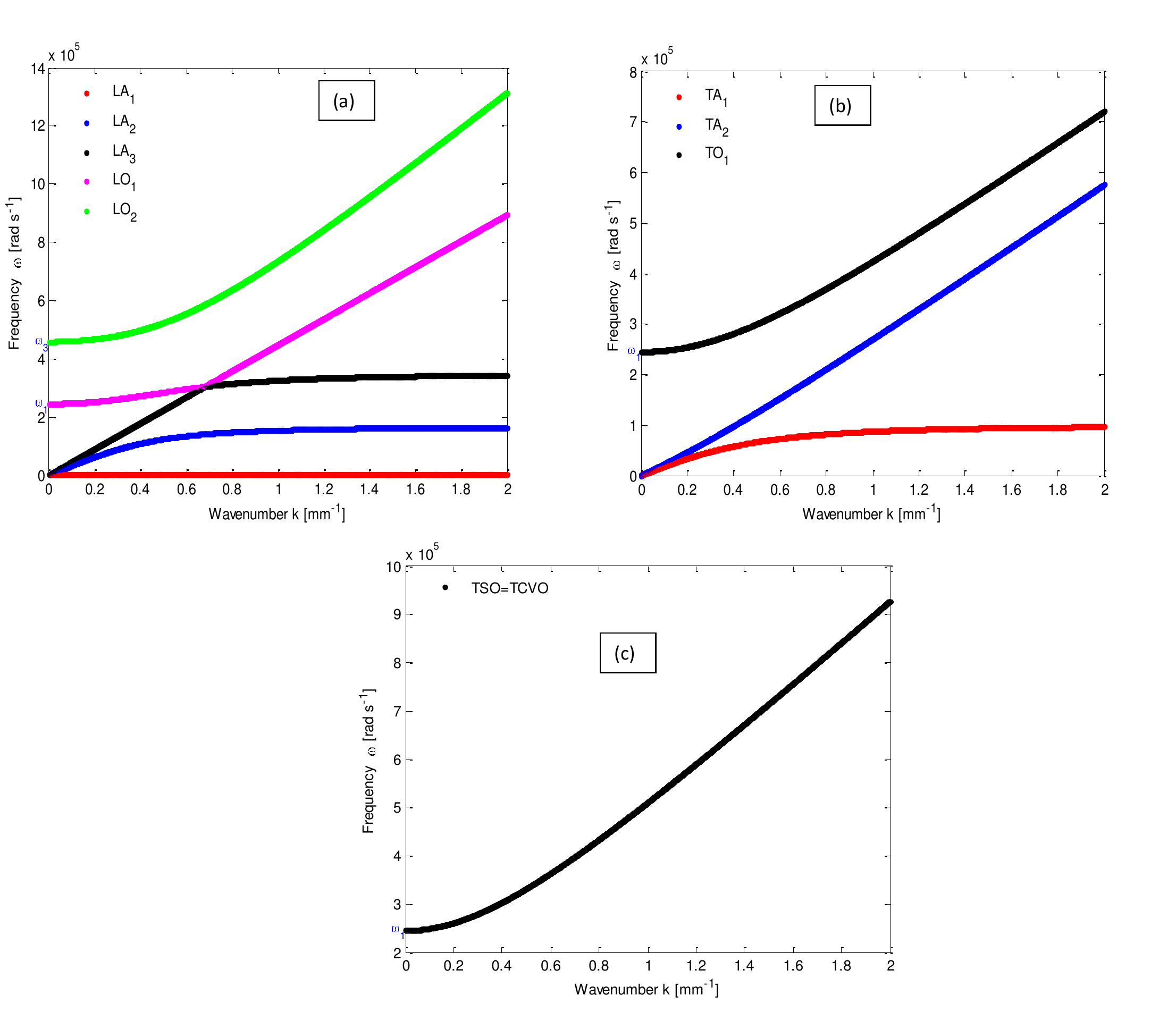}
			\caption{Dispersion curves in the thermoelastic relaxed micro-morphic model having $a_1>0, a_2=0,$ $a_3=\displaystyle \frac{a_1}{3}$ with zero Cosserat couple modulus $\mu_c = 0$. [(a) Longitudinal branches, (b) Transverse branches, (c) Uncoupled Optic branch]}
		\end{center}
	\end{figure*}
	\begin{figure*}
		\begin{center}
			\includegraphics[width=0.85\textwidth]{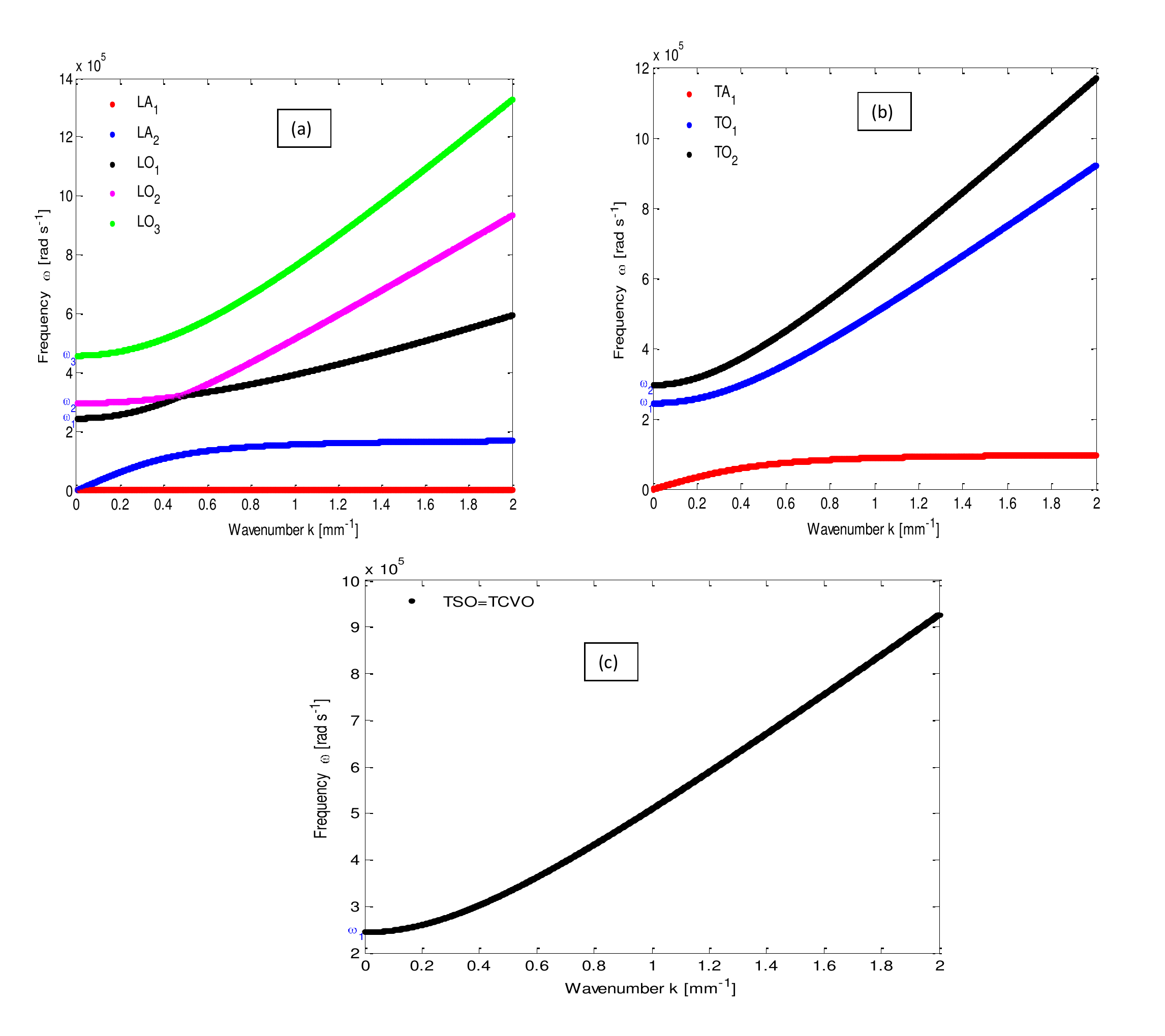}
			\caption{Dispersion curves in the  thermoelastic relaxed micro-morphic model having $a_1=a_2 >0,$ $a_3=0$ with positive Cosserat couple modulus $\mu_c > 0$. [(a) Longitudinal branches, (b) Transverse branches, (c) Uncoupled Optic branch]}
		\end{center}
	\end{figure*}
	Figures 6 and 7 depict the dispersion curves of various existing waves in Model-V and Model-VI, respectively. The effects of the Cosserat couple modulus $\mu_c$ on the dispersion curves can be noticed by comparing the respective graphs of Figures 6 and 7. On comparing these figures, one can notice that the presence of the Cosserat couple modulus $\mu_c$ is solely responsible for frequency band gap. On comparing Figure 6(a) with Figure 7(a), it is seen that there exist three-LO and two-LA branches for the  Cosserat couple modulus $\mu_c >0$, while there exist two-LO and three-LA branches for zero Cosserat couple modulus $\mu_c=0.$ A similar phenomena is seen while comparing Figure 6(b) with Figure 7(b). On comparing Figure 6(c) with Figure 7(c), we see that there are no effects of the  Cosserat couple modulus $\mu_c$ on TSO-TCVO branches.
	\begin{figure*}
		\begin{center}
			\includegraphics[width=0.85\textwidth]{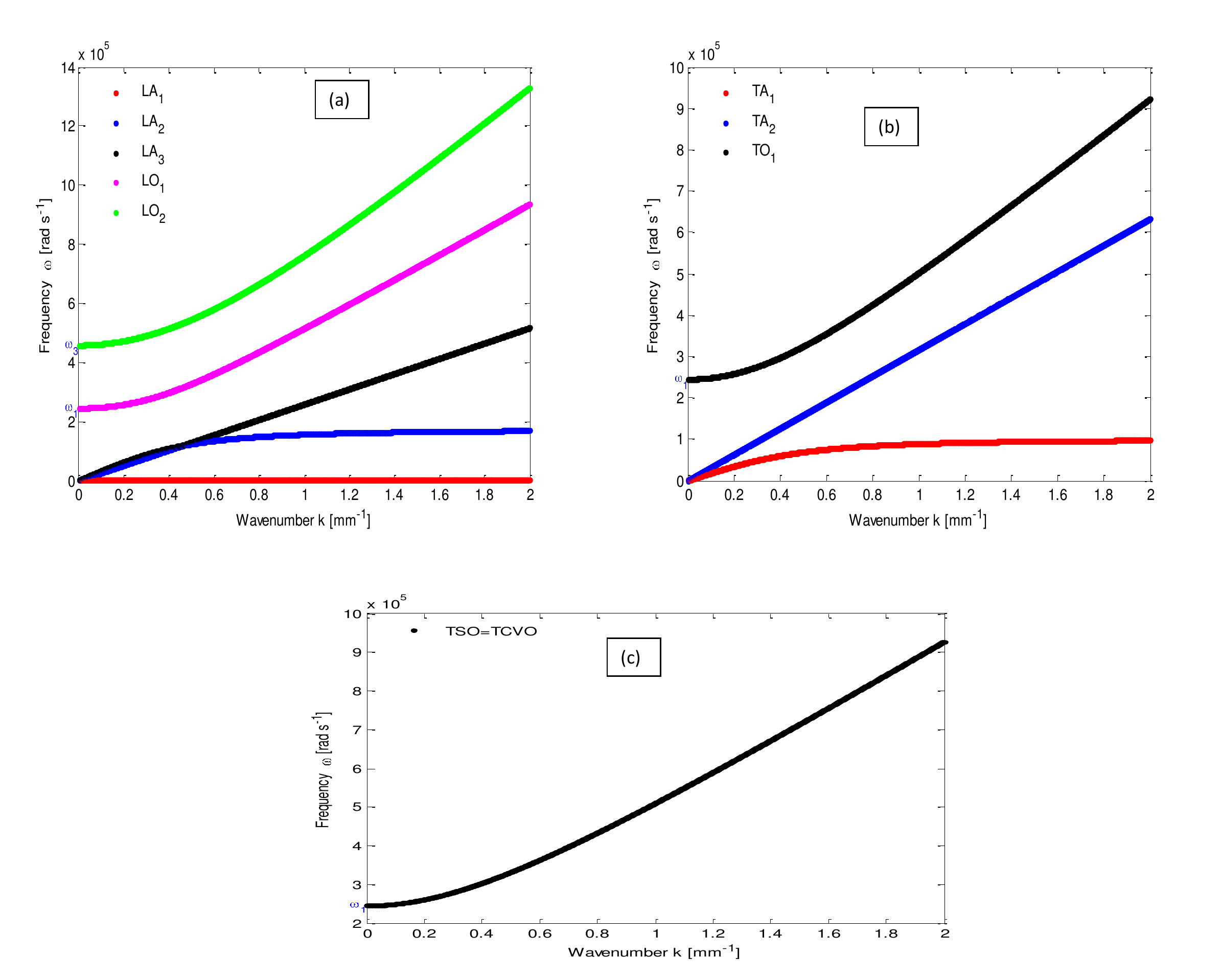}
			\caption{Dispersion curves in the thermoelastic relaxed micro-morphic model having $a_1=a_2 >0,$ $a_3=0$ with zero Cosserat couple modulus $\mu_c = 0$. [(a) Longitudinal branches, (b) Transverse branches, (c) Uncoupled Optic branches]}
		\end{center}
	\end{figure*}
	\section{Conclusions}
	The dispersion behavior of the various waves in the  thermoelastic relaxed micrmorphic medium has been analyzed under four different scenarios. It has been found that:
\begin{enumerate}
 	\item [(i)] the thermal effects give rise to new wave and generate coupling with longitudnal waves.
 		\item [(ii)] the new wave arising due to thermal effects is dispersive.
 	\item [(iii)]the models with positive Cosserat couple modulus model  $\mu_c >0$  contain band gaps and are found to be the same as in the   Madeo et al.\,\cite{madeo2015wave} model.
 	\item [(iv)] the models with zero Cossearat couple modulus $\mu_c =0$ contain no band gap.
 	\item [(v)] the transverse rotational acoustic wave of Madeo et al.\,\cite{madeo2015wave} transforms into a longitudinal wave due  to thermal effects. 
 	
 \end{enumerate}
\clearpage
 \footnotesize
\bibliography{bib_thermo,bib_agn,bib_main}

\end{document}